\documentclass[11pt]{article}
\usepackage{}
\usepackage[total={6in, 9in}]{geometry}
\usepackage{amsfonts}
\usepackage{amsthm}
\usepackage{amssymb}
\usepackage{amsmath,amsfonts}
\usepackage{mathrsfs}
\usepackage{cases}
\usepackage{latexsym,bm}
\usepackage{indentfirst}
\usepackage{color}
\usepackage{ifpdf}
\usepackage{graphicx}
\usepackage{psfrag}
\usepackage{dsfont}
\usepackage[
pdfauthor={CESS},
pdftitle={Spanning trails},
pdfstartview=XYZ,
bookmarks=true,
colorlinks=true,
linkcolor=blue,
urlcolor=blue,
citecolor=blue,
bookmarks=true,
linktocpage=true,
hyperindex=true
]{hyperref}

\newtheorem{THM}{\textbf{Theorem}}

\newtheorem{LEM}{\textbf{Lemma}}[section]
\newtheorem{CON}{\textbf{Conjecture}}

\newcommand{\pf}{\noindent\textbf{Proof}.\quad}

\newcommand{\iC}{\overset{\rightharpoonup }{C}}

\def\dist{{\fam0 dist}}

\begin{document}
\title{Spanning trails with maximum degree at most 4 in $2K_2$-free graphs}
\author{Guantao Chen\thanks{Department of Mathematics and
Statistics, Georgia State University, Atlanta, Georgia 30303, U.S.A.
        ({\tt gchen@gsu.edu}). }
	\and M. N. Ellingham\thanks{Department of Mathematics,
Vanderbilt University, Nashville, Tennessee 37240, U.S.A. ({\tt
mark.ellingham@vanderbilt.edu}).  Supported by Simons Foundation award
429625.}
        \and Akira Saito\thanks{Department of Information Science, Nihon University,
Tokyo 156--8550, JAPAN ({\tt asaito@chs.nihon-u.ac.jp}).
       }
\and Songling Shan\thanks {Department of Mathematics, Vanderbilt University, Nashville, Tennessee 37240,
U.S.A. ({\tt songling.shan@vanderbilt.edu}).}
}

\date{}
\maketitle

\emph{\textbf{Abstract}.}
A graph is called $2K_2$-free if it does not contain two independent edges as an induced subgraph. Mou and
Pasechnik conjectured that every $\frac{3}{2}$-tough $2K_2$-free graph with at least three vertices
has a spanning trail with maximum degree at most $4$. In this paper, we confirm this conjecture.
 We also provide examples for all $t < \frac{5}{4}$ of $t$-tough graphs
that do not have a spanning trail with maximum degree at most $4$.


\vspace{2mm}

\section{Introduction}

Graphs considered in this paper are simple, undirected,  and finite. Let $G$ be a graph.
 Let $V(G)$ and  $E(G)$ be the vertex set and edge set of $G$,
respectively. For $v\in V(G)$,  $N_G(v)$ denotes the set of neighbors
of $v$ in $G$, and $d_G(v)=|N_G(v)|$ the degree of $v$ in $G$. If
$S\subseteq V(G)$ then the subgraph induced by $V(G)-S$ is denoted by
$G-S$. For notational simplicity we write $G-\{x\}$ for $G-x$. Let
$u,v\in V(G)$ be two vertices. Then $\dist_G(u,v)$, the {\it distance
between $u$ and $v$ in $G$\/}, is defined to be the length of a shortest
path connecting $u$ and $v$ in $G$.
 If $uv\not\in E(G)$, we write $G+uv$ for the new graph obtained from
$G$ by adding the edge $uv$. If $uv\in E(G)$, then $G-uv$ denotes the
graph obtained from $G$ by deleting the edge $uv$. Let $V_1,
V_2\subseteq V(G)$ be two disjoint sets. Then $E_G(V_1,V_2)$ is the set
of edges of $G$ with one end in $V_1$ and the other end in $V_2$. The
graph $G$ is called $2K_2$-free if it does not contain two independent
edges as an induced subgraph.

The number of components of $G$ is denoted by $c(G)$. Let $t\ge 0$ be a
real number. The graph is said to be {\it $t$-tough\/} if $|S|\ge t\cdot
c(G-S)$ for each $S\subseteq V(G)$ with $c(G-S)\ge 2$. The {\it
toughness $\tau(G)$\/} is the largest real number $t$ for which $G$ is
$t$-tough, or is defined as $\infty$ if $G$ is complete. This concept, a
measure of graph connectivity and ``resilience'' under removal of
vertices, was introduced by Chv\'atal~\cite{chvatal-tough-c} in 1973.
 It is  easy to see that  if $G$ has a hamiltonian cycle
then $G$ is 1-tough. Conversely,
 Chv\'atal~\cite{chvatal-tough-c}
 conjectured that
there exists a constant $t_0$ such that every
$t_0$-tough graph is hamiltonian.
 Bauer, Broersma and Veldman~\cite{Tough-counterE} have constructed
$t$-tough graphs that are not hamiltonian for all $t < \frac{9}{4}$, so
$t_0$ must be at least $\frac{9}{4}$.

There are a number of papers on
 Chv\'atal's toughness conjecture,
and it  has
been verified when restricted to a number of graph
classes~\cite{Bauer2006},
including planar graphs, claw-free graphs, co-comparability graphs, and
chordal graphs.
Recently, Broersma, Patel and Pyatkin~\cite{2k2-tough} proved that
every 25-tough $2K_2$-free graph on at least three vertices is
hamiltonian.

Another direction inspired by  Chv\'atal's toughness conjecture
is investigating the existence of spanning  substructures weaker than
hamiltonian cycles
for a given toughness. For example,
$k$-trees,   $k$-walks, and $k$-trails are substructures of this kind.
 Let $k$ be a positive integer.
 A {\it $k$-tree\/} is a tree with maximum degree at most $k$, and a {\it
$k$-walk\/} is a closed walk with each vertex repeated at most $k$ times.
 A $k$-walk can be obtained from a $k$-tree by visiting each
edge of the tree twice.
 A {\it $k$-trail\/} is a $k$-walk
with no repetition of edges.  A graph has a spanning $k$-trail if and only if
it has a spanning Eulerian subgraph with maximum degree at most $2k$.
 A spanning 2-tree is just a hamiltonian path and a spanning
1-walk/1-trail is a hamiltonian cycle.

 In 1990, Jackson and
Wormald~\cite{JW-k-walks} made the following conjecture.

\begin{CON}\label{k-walkc}
Let $k\ge 2$ be a positive integer. Then every $\frac{1}{k-1}$-tough graph has a spanning $k$-walk.
\end{CON}

 Mou and Pasechnik~\cite{1412.0514,2k21} confirmed Jackson and Wormald's conjecture for
$2K_2$-free graphs. In~\cite{2k21}, they proposed the following two conjectures.

\begin{CON}\label{2k2t}
Every $\frac{3}{2}$-tough $2K_2$-free graph with at least three vertices has a spanning 2-trail.
\end{CON}

\begin{CON}\label{2k2h}
Every $2$-tough $2K_2$-free graph with at least three vertices is hamiltonian.
\end{CON}

The class of $2K_2$-free graphs is well studied, for instance, see~\cite{
2k2-tough, CHUNG1990129, MR845138, MR2279069, 1412.0514, 2k21,MR1172684}.
It is a superclass of {\it split\/} graphs,
where the vertices can be partitioned into a clique and an independent set.
One can also easily check that every {\it cochordal\/} graph (i.e., a graph that is the complement of a
chordal graph) is $2K_2$-free and so the class of $2K_2$-free graphs is
at least as rich as the class of chordal graphs.

In this paper, we confirm Conjecture~\ref{2k2t}.

\begin{THM}\label{main}
Let $G$ be a $\frac{3}{2}$-tough $2K_2$-free graph with at least three vertices. Then $G$
has a spanning 2-trail.
\end{THM}

There is a large literature proving the existence of a spanning closed
trail under various conditions; a graph with a spanning closed trail is
called {\it supereulerian}.  A recent paper in this area, providing
references to other papers, is \cite{MR3442576}.
 However, apart from results on hamiltonicity there do not seem to be
many results on spanning closed trails with bounded degree.  Other than
Theorem \ref{main}, the only one we are aware of is in \cite{MR2040107},
which proves that a $2$-edge-connected $n$-vertex graph $G$ with $n \ge
7$ and $\sigma_3(G) \ge n$ has a spanning $2$-trail, where $\sigma_3(G)$
is the minimum degree sum over all triples of pairwise independent
vertices.

We prove Theorem~\ref{main} in Section 2. In Section 3, we construct $2K_2$-free
graphs with toughness close to $\frac{5}{4}$ but containing no spanning 2-trail.

\section{Proof of Theorem~\ref{main}}

We need the following lemma in proving Theorem~\ref{main}.

\begin{LEM}\label{tough}
Let $G$ be a bipartite graph with  partite sets $X$ and $Y$.
If for every $S\subseteq X$, $|N_G(S)|\ge \frac{3}{2}|S|$,
then $G$ has a subgraph $H$  covering $X$\,(meaning that $X\subseteq V(H)$)
such that for every $x\in X$, $d_H(x)=2$ and for every $y\in Y$,
$d_H(y)\le 2$.
\end{LEM}

\pf Form $G'$ from $G$ by replacing each $x\in X$ by $x_1, x_2, x_3$,
each $y\in Y$ by $y_1,y_2$, and each $xy\in E(G)$ by six edges $x_iy_j$,
$1\le i\le 3$, $1\le j\le 2$. Let $\pi$ be the natural projection from
$G'$ to $G$ with $$ \pi(x_i)=x,\quad \pi(y_j)=y, \quad \pi(x_iy_j)=xy.
$$
 Let $X'=\pi^{-1}(X)$ be the inverse image of $X$ under $\pi$.
 For each $S'\subseteq X'$, let $S=\pi (S')$.
Then $|S'|\le 3|S|\le 2|N_G(S)|=|N_{G'}(S')|$.
 Thus, by Hall's Theorem, $G'$ has a matching $M'$ covering $X'$.
 The projection  $\pi(M')$ of $M'$ is a graph containing all the
vertices in $X$ such that each vertex in $X$ has degree 2 or 3, and each
vertex in $Y$ has degree at most 2.
 In  $\pi(M')$, for each $x\in X$ with degree 3, delete one edge
incident to $x$. Then the graph $H$ induced by the remaining edges
is the desired graph.
 \qed

 \smallskip
 We cannot reduce the number $\frac32$ in Lemma \ref{tough}.
 To see this, take $k \ge 1$, $X$ with $|X|=2k$, and $Y
= Y_1 \cup Y_2$ with $|Y_1| = 2k$ and $|Y_2|=k$.  To form $G$, join each
vertex of $X$ to a distinct vertex of $Y_1$ (giving a matching) and join
every vertex of $X$ to every vertex of $Y_2$.  Then $G$ has a subgraph
$H$ as described, but if we delete any $y \in Y$ then no such
subgraph exists although $G-y$ satisfies the condition of Lemma \ref{tough}
with $\frac{3k-1}{2k}$ instead of $\frac{3}{2}$.

A subgraph $G^*\subseteq G$ is called {\it dominating\/}
if $G-V(G^*)$ is an edgeless graph. Mou and Pasechnik
proved the existence of a dominating cycle in $2K_2$-free graphs.
In fact, the proof of \cite[Theorem 3]{1412.0514} implies the following.

\begin{LEM}\label{dominating}
Let $G$ be a $2K_2$-free graph containing a cycle.
Then some longest cycle of $G$ is dominating.
\end{LEM}

 \noindent
 {\bf Proof of Theorem~\ref{main}.}
 As $G$ is $\frac{3}{2}$-tough, $G$
is 3-connected. So $G$ has a cycle. Let $C$ be a dominating longest
cycle of $G$, which exists by Lemma~\ref{dominating}.
 Let $\iC$ denote a forward orientation of $C$. For
a vertex $x\in V(C)$, we let $x^+$ denote the successor of $x$ on $\iC$,
and if $S \subseteq V(C)$ we define $S^+ = \{x^+ \;|\; x \in S\}$.
We may assume $V(G)-V(C)\ne \emptyset$. Otherwise, $C$ is a spanning
1-trail.

 \smallskip
 \noindent
{\bf  Claim A:} Let $x\in V(G)-V(C)$.

 \vspace{-\parskip}
 \vspace{-\topsep}
 \begin{itemize}
 \setlength{\itemsep}{0pt}
  \item [(a)] $N_G(x)$ does not contain two consecutive vertices on $C$.
  \item [(b)] If $y, z \in N_G(x)$ with $y \ne z$ then there is no  path
from $y^+$ to $z^+$ that is that is internally disjoint from $C$; in
particular, $y^+z^+ \notin E(G)$.
  \item [(c)] $C$ has at least $7$ vertices.
 \end{itemize}
 \vspace{-\topsep}

\pf Both (a) and (b) follow by standard arguments.
 We only prove (c) here.
 Since $G-V(C)$ is edgeless, $N_G(x) \subseteq V(C)$.
 By (a), $N_G(x)^+$ is disjoint from $N_G(x)$.
 As $G$ is $\frac{3}{2}$-tough, $\delta(G)\ge 3$, so $|N_G(x)| =
|N_G(x)^+| \ge 3$.  Thus, $|V(C)| \ge 6$, and $|V(C)|=6$ precisely when
$|N_G(x)|=3$ and $V(C) = N_G(x) \cup N_G(x)^+$.
 In that case, by (b) the vertices of $N_G(x)^+$ belong to separate
components in $G-N_G(x)$.
 Thus, $c(G-N_G(x))\ge 4$, and so $\frac{|N_G(x)|}{c(G-N_G(x))}\le
\frac{3}{4}<\frac{3}{2}$, contradicting the toughness of $G$. \qed

 \smallskip
Let $G'=G-E(G[V(C)])$ with partite sets $X=V(G)-V(C)$ and $Y=V(C)$.
Since $G$ is $\frac{3}{2}$-tough and $X$ is an independent set in $G$,
we have that for any $S\subseteq X$, $|N_{G'}(S)|\ge \frac{3}{2}|S|$
 (even when $|S|=1$, because then $c(G-N_{G'}(S)) \ge 2$ by (a) of Claim
A).
 Applying  Lemma~\ref{tough} to $G'$, we see that $G'$ (hence $G$) has a
subgraph $H$ such that for any $x\in X$, $d_H(x)=2$ and for any $y\in
Y\cap V(H)$, $d_H(y)=1$ or $d_H(y)=2$. Subject to this property, we
choose a subgraph $H$ of $G$ such that the number of components  in $H$
is smallest.
 Let $H_1,\cdots, H_\ell$ be the components of $H$. Each $H_i$ is either
a path or a cycle. Assume, without loss of generality, that $H_1,\cdots,
H_{p}$ are paths and $H_{p+1}, \cdots, H_\ell$ are cycles. For each path
$H_i$\,($1\le i\le p$), let $u_i$ and $v_i$ denote its endvertices
(these two vertices are on $C$ by the construction of $H$). Let $s_i$
and $t_i$ denote the neighbor of $u_i$ and $v_i$ in $H$, respectively.
Note that $s_i$ and $t_i$ are vertices from $V(G)-V(C)$ and $s_i=t_i$ if
$H_i$ has length 2. Note also that $C\cup\left(\bigcup_{p+1 \le i \le
\ell} H_i\right)$  is a spanning  2-trail if $p=0$. Therefore, we
assume $p\ge 1$.

 \smallskip
 \noindent
{\bf  Claim B:} Each of the following holds.

 \vspace{-\parskip}
 \vspace{-\topsep}
 \begin{itemize}
 \setlength{\itemsep}{0pt}%
  \item [(a)] $s_iu_j, s_iv_j,  t_iu_j, t_iv_j \not\in E(G)$, for all
$i,j$ with $i\ne j$ and $i,j\in \{1,\cdots, p\}$.
  \item [(b)] Let $u$ be an endvertex of $H_i$ and $v$ be an endvertex
of $H_j$, where $i\ne j$ and $i,j\in \{1,\cdots, p\}$.   Then $uv\in
E(G)$.
 \end{itemize}
 \vspace{-\topsep}

\pf For (a), if say $s_i u_j \in E(G)$ then we could replace $s_iu_i$
by $s_iu_j$ in $H$ to obtain fewer components.
 For (b), let $s$ be the neighbor of
$u$ on $H_i$, and $t$ be the neighbor of $v$ on $H_j$.
Note that $s,t\in V(G)-V(C)$.
Since $i\ne j$,
we have $s\ne t$. By (a),
we have $sv, tu\not\in E(G)$.
Furthermore, $st\not\in E(G)$ as
$G-V(C)$
is edgeless. So $uv\in E(G)$ by the
$2K_2$-freeness
of $G$.
\qed

 \smallskip
 \noindent
{\bf  Claim C:} Let $q$ be an integer with $1\le q\le p$,
and let $V_q=\bigcup_{1\le i\le q}V(H_i)$.
Then $G[V_q]-E(C)$ contains a path $P_q$ with vertex set $V_q$
such that for each $i$ with $1\le i\le q$, $H_i$ is a subpath of $P_q$
and both endvertices of  $P_q$ belong to $\{u_1,\cdots, u_q, v_1,\cdots, v_q\}$.

\pf
We show this claim by induction on $q$.
For $q=1$, $H_1$ itself is a desired path. So we assume that $q\ge 2$.
By the induction hypothesis,
$G[V_{q-1}]-E(C)$ contains a path $P_{q-1}$ with the desired property.
Assume, without loss of generality, that the two endvertices of $P_{q-1}$ are $u_a$ and $v_b$
with $a,b \in \{1,\cdots, q-1\}$.
As $|V(C)|\ge 7$ by (b) of Claim A, we see that one of $\dist_C(u_{a},u_{q}), \dist_C(u_{a},v_{q}),
\dist_C(v_{b},u_{q}), \dist_C(v_{b},v_{q})$ must be at least 2.
Assume, without loss of generality,   that $\dist_C(u_{a},v_{q})\ge 2$.
Then $u_av_q\in E(G)$ by (b) of Claim B and
 $u_av_q\in E(G)-E(C)$ since $\dist_C(u_{a},v_{q})\ge 2$.
Thus, $P_{q-1}\cup H_q+u_av_q$ is a desired path.
\qed

 Let $D = \bigcup_{p+1 \le i \le \ell} H_i$ be the union of the cycle
components of $H$.  Consider two cases.

 \smallskip\noindent
{\bf Case 1:} $p\ge 2$.

Let $P_p$ be a path with the property stated in Claim C.
Assume, without loss of generality, that the endvertices of $P_p$
are $u_1$ and $v_p$.
By (b) of Claim B,  we have $v_pu_1\in E(G)$.
Let
$$
T=\left\{
    \begin{array}{ll}
      C\cup D \cup P_p-v_pu_1, & \hbox{if $v_pu_1\in E(C)$;} \\
      C\cup D \cup P_p+v_pu_1, & \hbox{if $v_pu_1\in E(G)-E(C)$.}
    \end{array}
  \right.
$$
Then $T$ is a spanning 2-trail of $G$.

 \smallskip\noindent
{\bf Case 2:} $p= 1$.

Assume first that $|V(H_1)|\ge 4$. Consider the two edges $s_1u_1$ and $t_1v_1$.
Again,  we have $\{s_1v_1, t_1u_1, u_1v_1\}\cap  E(G)\ne \emptyset$
by the $2K_2$-freeness of $G$. Let
$$
T=\left\{
    \begin{array}{ll}
      C\cup H-s_1u_1+s_1v_1, & \hbox{if $s_1v_1\in E(G)$;} \\
      C\cup H-t_1v_1+t_1u_1, & \hbox{if $t_1u_1\in E(G)$;} \\
      C\cup H-u_1v_1, & \hbox{if $u_1v_1\in E(C)$;} \\
      C\cup H+u_1v_1, & \hbox{if $u_1v_1\in E(G)-E(C)$.}
    \end{array}
  \right.
$$
Then $T$ is a spanning 2-trail of $G$.

 Assume now that $|V(H_1)|=3$.  Suppose that $\dist_C(u_1, v_1)
\ge 3$.
 As $u_1^+v_1^+\not\in E(G)$ by (b) of Claim A, we
have $\{u_1v_1, u_1v_1^+, v_1u_1^+\}\cap E(G)\ne \emptyset$ by the
$2K_2$-freeness of $G$. Note that $\{u_1v_1, u_1v_1^+, v_1u_1^+\}\cap
E(C)=\emptyset$  as $\dist_C(u_1, v_1)\ge 3$. Let $$ T=\left\{
    \begin{array}{ll}
      C\cup H+u_1v_1^+-v_1v_1^+, & \hbox{if $u_1v_1^+\in E(G)$;} \\
      C\cup H+v_1u_1^+-u_1u_1^+, & \hbox{if $v_1u_1^+\in E(G)$;}\\
      C\cup H+u_1v_1, & \hbox{if $u_1v_1\in E(G)$.}
    \end{array}
  \right.
$$
 In the first case the vertex  $v_1^+$  may also be contained in
$D$, but when we add the edge $u_1v_1^+$
and remove the edge $v_1v_1^+$, the degree of $v_1^+$ in $T$
is the same as in $C\cup H$. The same applies to $u_1^+$ in the
second case.
 Thus the degree of each vertex in $T$ is at most 4, and  $T$ is a
spanning 2-trail of $G$.

 Suppose that $N_G(s_1) - V(H) \ne \emptyset$.  Then
$N_G(s_1)-V(D)$, which includes $u_1$ and $v_1$, contains at least three
vertices.
 By Claim A, these vertices are pairwise nonadjacent and $|V(C)| \ge 7$,
so there are $u', v' \in N_G(s_1)-V(D)$ with $\dist_C(u', v') \ge 3$.
We replace $H_1$ by the path $u' s_1 v'$ and apply the argument above.

 Therefore, we assume all neighbors of $s_1$ not in $H_1$ lie in $D$,
which must be nonempty.
 Suppose that $N_G(x') \subseteq V(H)$ for all $x' \in X \cap V(D)$.
 Then deleting all the $|X|+1$ neighbors of vertices in $X$ on $C$
results in at least $|X|$ components.
 Since $|V(D) \cap X| \ge 2$ and $s_1 \in X$, $|X| \ge 3$, so
 $\frac{|X|+1}{|X|}\le \frac{4}{3}<\frac{3}{2}$,
 contradicting the toughness of $G$.
 Therefore there exist $x'$ and $u'$ with $x'\in X \cap V(D)$ and $u'
\in  N_G(x')-V(H)$.
 Let $x'v' \in E(D)$ and $D'=D-x'v'+x'u'$. Replacing $D$ by $D'$ in
$H$, we see that the new graph has the same property as $H$, but it has
two components that are paths, so we may apply Case 1.

 The proof of Theorem \ref{main} is now complete.
\qed

\section{An Extremal Example}

In this section, we construct a family of $2K_2$-free graphs
with toughness approaching $\frac{5}{4}$
that do not contain any spanning 2-trail.

Let $n\ge 2$ be an integer, $Q_1=K_{4n}$,
the complete graph on $4n$ vertices, $Q_2=\overline{K_{4n}}$,
the empty graph on $4n$ vertices, and $Q_3=K_{n-1}$. Let
$G_n$ be a graph with $V(G_n)=V(Q_1)\cup V(Q_2)\cup V(Q_3)$ and
$E(G_n)$ consisting of all edges in $Q_1$ and $Q_3$, all edges between $V(Q_3)$
and $V(Q_1)\cup V(Q_2)$, and a perfect matching between $Q_1$ and $Q_2$.
It is easy to check that $G$ is $2K_2$-free.

We claim that $\lim\limits_{n\rightarrow\infty}\tau(G_n)=\frac{5}{4}$.
Let $S\subseteq V(G_n)$ be a cutset such that $\tau(G_n)=\frac{|S|}{c(G_n-S)}$.
Then $Q_3\subseteq S$ as each vertex in $Q_3$ is adjacent to every other
vertex of $G_n$. Also, $S\cap V(Q_2)=\emptyset$.
Otherwise, as $c(G-(S-V(Q_2))) \ge c(G-S)$, we get $\frac{|S-V(Q_2)|}{c(G_n-(S-V(Q_2)))}<\frac{|S|}{c(G_n-S)}=\tau(G_n)$,
contradicting the toughness of $G$.
Thus, $c(G-S)=|S\cap V(Q_1)|+1$ if $V(Q_1)\not\subseteq S$ and $c(G-S)=4n$ otherwise.
In the latter case, $\frac{|S|}{c(G_n-S)}=\frac{5n-1}{4n}$.
So assume $V(Q_1)\not\subseteq S$ and $|V(Q_1)\cap  S|=r$, where $1\le r\le 4n-1$.
Then $\frac{n-1+r}{r+1}$ is a decreasing function of $r$ which achieves its minimum when
$r=4n-1$. Hence,  $\tau(G_n)=\frac{|S|}{c(G_n-S)}=\frac{5n-2}{4n}$, which approaches
$\frac{5}{4}$  as $n \to \infty$.

We show now that $G_n$ has no spanning 2-trail. Suppose on the contrary that
$T$ is a spanning 2-trail of $G_n$. Let $v\in V(Q_2)$ be a vertex. Then $d_T(v)\ge 2$.
As $|N_G(v)\cap V(Q_1)|=1$, $|N_T(v)\cap V(Q_3)|\ge 1$.
Thus, $|E_T(V(Q_3), V(Q_2))|\ge 4n$.
Since $|V(Q_3)|=n-1$, by the Pigeonhole Principle
there is a vertex from $Q_3$ that has degree at least 5 in $T$.
This contradicts the assumption that $T$ is a 2-trail.

From the example above, we suspect the following might be true.

\begin{CON}
Any $\frac{5}{4}$-tough $2K_2$-free graph with at least three vertices has a spanning 2-trail.
\end{CON}

 Our proof of Theorem \ref{main} relies on Lemma \ref{tough}, which
cannot be improved, so a new strategy will be needed to obtain a
positive answer to this conjecture.


\bibliographystyle{plain}
\bibliography{SSL-BIB}

\begin{thebibliography}{10}

\bibitem{Tough-counterE}
D.~Bauer, H.~J. Broersma, and H.~J. Veldman.
\newblock Not every 2-tough graph is {H}amiltonian.
\newblock In {\em Proceedings of the 5th {T}wente {W}orkshop on {G}raphs and
  {C}ombinatorial {O}ptimization ({E}nschede, 1997)}, volume~99, pages
  317--321, 2000.

\bibitem{Bauer2006}
D.~Bauer, H.J. Broersma, and E.~Schmeichel.
\newblock Toughness in graphs -- a survey.
\newblock {\em Graphs and Combinatorics}, 22(1):1--35, 2006.

\bibitem{2k2-tough}
H.~Broersma, V.~Patel, and A.~Pyatkin.
\newblock On toughness and {H}amiltonicity of {$2K_2$}-free graphs.
\newblock {\em J. Graph Theory}, 75(3):244--255, 2014.

\bibitem{CHUNG1990129}
F.~R.~K. Chung, A.~Gy{\'a}rf{\'a}s, Z.~Tuza, and W.~T. Trotter.
\newblock The maximum number of edges in {$2K_2$}-free graphs of bounded
  degree.
\newblock {\em Discrete Math.}, 81(2):129--135, 1990.

\bibitem{chvatal-tough-c}
V.~Chv{\'a}tal.
\newblock Tough graphs and {H}amiltonian circuits.
\newblock {\em Discrete Math.}, 5:215--228, 1973.

\bibitem{MR845138}
M.~El-Zahar and P.~Erd{\H{o}}s.
\newblock On the existence of two nonneighboring subgraphs in a graph.
\newblock {\em Combinatorica}, 5(4):295--300, 1985.

\bibitem{MR2040107}
M.~N. Ellingham, X.~Zha, and Y.~Zhang.
\newblock Spanning 2-trails from degree sum conditions.
\newblock {\em J. Graph Theory}, 45(4):298--319, 2004.

\bibitem{JW-k-walks}
B.~Jackson and N.~C. Wormald.
\newblock {$k$}-walks of graphs.
\newblock {\em Australas. J. Combin.}, 2:135--146, 1990.
\newblock Combinatorial mathematics and combinatorial computing, Vol. 2
  (Brisbane, 1989).

\bibitem{MR3442576}
Ping Li, Hao Li, Ye~Chen, Herbert Fleischner, and Hong-Jian Lai.
\newblock Supereulerian graphs with width {$s$} and {$s$}-collapsible graphs.
\newblock {\em Discrete Appl. Math.}, 200:79--94, 2016.

\bibitem{MR2279069}
D.~Meister.
\newblock Two characterisations of minimal triangulations of {$2K_2$}-free
  graphs.
\newblock {\em Discrete Math.}, 306(24):3327--3333, 2006.

\bibitem{2k21}
G.~Mou and D.~Pasechnik.
\newblock On $k$-walks in $2{K}_2$-free graphs.
\newblock {\em arXiv:1412.0514v2}, 2014.

\bibitem{1412.0514}
G.~Mou and D.~Pasechnik.
\newblock Edge-dominating cycles, k-walks and hamilton prisms in $2{K}_2$-free
  graphs.
\newblock {\em arXiv:1412.0514v4}, 2015.

\bibitem{MR1172684}
M.~Paoli, G.~W. Peck, W.~T. Trotter, Jr., and D.~B. West.
\newblock Large regular graphs with no induced {$2K_2$}.
\newblock {\em Graphs Combin.}, 8(2):165--197, 1992.

\end{thebibliography}

\end{document}